\documentclass{article}

\usepackage{graphicx} 
\usepackage{amsmath,amssymb}
\usepackage{amsthm}
\usepackage{authblk}

\usepackage{url} 
\usepackage[breaklinks,colorlinks]{hyperref}

\theoremstyle{remark}
\newtheorem{remark}{Remark}

\usepackage{listings}
\newcommand{\func}[1]{\texttt{#1}}
\newcommand{\var}[1]{\texttt{#1}}
\newcommand{\const}[1]{\texttt{#1}}
\newcommand{\code}[1]{\texttt{#1}}

\title{An algorithm for scaling vectors by the reciprocal of a complex number}
\author[]{Weslley da Silva Pereira}
\affil[]{Department of Mathematical and Statistical Sciences, University of Colorado Denver\\ \href{mailto://weslley.pereira@ucdenver.edu}{weslley.pereira@ucdenver.edu}}

\begin{document}

\maketitle

\begin{abstract}
    This document describes an algorithm to scale a complex vector by the reciprocal of a complex value. The algorithm computes the reciprocal of the complex value and then scales the vector by the reciprocal. Some scaling may be necessary due to this 2-step strategy, and the proposed algorithm takes scaling into account. This algorithm is supposed to be faster than the naive approach of dividing each entry of the vector by the complex value, without losing much accuracy. It also serves as a single strategy for scaling vectors by the reciprocal of a complex value, which improves the software maintainability.
\end{abstract}

\section{Introduction}

In LAPACK, the routines \func{\{S,D,CS,ZD\}RSCL} scale a complex vector $x$ by the reciprocal of a real number $a$. This is equivalent to the following code:
\begin{lstlisting}[language=Fortran]
do i = 1, n
    x(i) = x(i) / a
end do
\end{lstlisting}
More specifically, the routines \func{RSCL} do two extra things:
\begin{enumerate}
\item They use the BLAS routine \func{SCAL} to do the scaling of the reciprocal of $a$. Thus, they use only one division and $n$ multiplications instead of the $n$ divisions in the code above.
\item The operation $(1/a) \cdot x_i$ can overflow or underflow in finite precision, even when $x_i/a$ does not. Therefore, the routines \func{RSCL} also check the range of $a$ and scale $x$ by a power of two if necessary. Thus, depending on the value of $a$, the code will also do extra $n$ multiplications by a power of two.
\end{enumerate}

As a result, \func{RSCL} are usually faster than the code above. The simple code above, however, is usually more accurate than \func{RSCL} since it uses fewer floating-point operations. When $x$ is complex, the performance gain is notable since complex divisions are always more expensive than complex multiplications, if we exclude the trivial case where the denominator has zero real or imaginary parts. Moreover, the accuracy loss in the complex case can be negligible since computing a complex division naturally translates to computing a complex multiplication of the numerator by the reciprocal of the denominator.

Some LAPACK routines, like \func{CLARFG} and \func{CGETF2}, need a reciprocal scaling where both $a$ and $x$ are complex. Since there is no LAPACK routine that does this, \func{CLARFG} and \func{CGETF2} have their own (and distinct) way to treat the reciprocal scaling. In fact, one finds the following code in \func{CLARFG} in LAPACK 3.11.0:
\begin{lstlisting}[language=Fortran]
ALPHA = CLADIV( CMPLX( ONE ), ALPHA-BETA )
CALL CSCAL( N-1, ALPHA, X, INCX )
\end{lstlisting}
where \func{CLADIV} is a LAPACK routine that computes the division of two complex numbers. In LAPACK 3.11.0, \func{CGETF2} has the following code:
\begin{lstlisting}[language=Fortran]
      IF( ABS(A( J, J )) .GE. SFMIN ) THEN
         CALL CSCAL( M-J, ONE / A( J, J ), A( J+1, J ), 1 )
      ELSE
         DO 20 I = 1, M-J
            A( J+I, J ) = A( J+I, J ) / A( J, J )
   20    CONTINUE
      END IF
\end{lstlisting}
where \const{SFMIN} is the smallest normal number such that $1/\const{SFMIN}$ does not overflow. Both codes perform the same operation, although the outcome can be different. For instance, if the absolute value of $\var{ALPHA-BETA}$ is very small in \func{CLARFG}, \func{CLADIV( CMPLX( ONE ), ALPHA-BETA )} may overflow. The overflow does not happen in \func{CGETF2}. Moreover, \func{CGETF2} may use $n$ complex divisions while \func{CLARFG} always uses a single one. These and other inconsistencies could be solved if there was a LAPACK routine that scales a vector by the reciprocal of a complex factor $a$.

\section{The algorithm}

This section presents an algorithm for the scaling of a complex vector \texttt{X} by the reciprocal of a complex number \texttt{A}. Let \texttt{AR} and \texttt{AI} be the real and imaginary parts of \texttt{A}, respectively. The algorithm is as follows:
\begin{enumerate}
  \item If \var{AI} is zero, then call \func{RSCL} (currently in LAPACK) using \var{AR}.
  \item If \var{AR} is zero, then if \var{AI} is in the safe range, call \func{SCAL} with the number \code{CMPLX( ZERO, -ONE / AI )}. Otherwise, do the proper scaling by a power of two in a way similar to what is done in \func{CSRSCL}.
  \item If both the real and imaginary parts of \var{A} are nonzero, then compute \var{UR} and \var{UI} such that \code{CMPLX( ONE / UR, -ONE / UI )} is the reciprocal of \var{A}, i.e.,
  \begin{lstlisting}[language=Fortran]
UR = AR + AI * ( AI / AR )
UI = AI + AR * ( AR / AI )
  \end{lstlisting}
  Note that \var{UR} and \var{UI} are always different from zero. NaNs only appear if either:
  \begin{itemize}
    \item \var{AR} or \var{AI} is a NaN.
    \item \var{AR} and \var{AI} are both infinite, in which case it makes sense to propagate a NaN.
  \end{itemize}
  We end up with three cases:
  \begin{enumerate}
    \item If \var{UR} and \var{UI} are both in the safe range, then call \func{SCAL} with the number \code{CMPLX( ONE / UR, -ONE / UI )}.
    \item If the absolute values of either \var{UR} or \var{UI} are smaller than \const{SFMIN}, it means that both \var{UR} and \var{UI} are small numbers (see Remark~\ref{rem:smallUR}). Therefore, it makes sense to scale \var{X} by \code{CMPLX( SFMIN / UR, - SFMIN / UI )}, and then scale the result by \code{ONE / SFMIN}. Use \func{SCAL} to do both scalings.
    \item If the absolute values of either \var{UR} or \var{UI} are greater than \code{ONE / SFMIN}, then check for infinite values:
    \begin{enumerate}
      \item If either \var{AR} or \var{AI} is infinite, then \var{UR} and \var{UI} are either both infinite or both NaNs. Therefore, scale \var{X} by \code{CMPLX( ONE / UR, -ONE / UI )} using \func{SCAL}.
      \item If either \var{UR} or \var{UI} is infinite and \var{AR} and \var{AI} are finite, then the algorithm generated infinite numbers that can be avoided. In this case, recompute scaled versions of \var{UR} and \var{UI} using \const{SFMIN}. Then, scale \var{X} by \const{SFMIN}, and then scale the result by \code{CMPLX( ONE / UR, -ONE / UI )}. Use \func{SCAL} to do both scalings.
      \item If \var{AR}, \var{AI}, \var{UR} and \var{UI} are finite, scale \var{X} by \const{SFMIN}, and then scale the result by \code{CMPLX( ONE / (SFMIN*UR), -ONE / (SFMIN*UI) )}. Use \func{SCAL} to do both scalings.
    \end{enumerate}
  \end{enumerate}
\end{enumerate}

Some comments about the algorithm:
\begin{itemize}
  \item It does not use complex divisions. Instead, it uses $n$ complex multiplications, $n$ real multiplications, and a fixed number of real multiplications and divisions.
  \item It propagates NaNs in \var{AR} and \var{AI} to the result. Moreover, it propagates NaNs if \var{AR} and \var{AI} are both infinite. It does not generate NaNs in other cases.
  \item It does not generate infinity numbers unless \var{A} is zero.
  \item It avoids flushing to zero in both the real and imaginary parts of the reciprocal of \var{A}. This is done by treating the special cases where \var{UR} or \var{UI} are big or infinity numbers.
  \item It does not lose accuracy due to subnormal numbers when generating the reciprocal of \var{A}.
\end{itemize}

\begin{remark}
\label{rem:smallUR}
    Assume \var{ABS(UR)} is smaller than \const{SFMIN} and both \var{AR} and \var{AI} are nonzero. Notice that \var{ABS(UR) = ABS(AR) + AI * ( AI / ABS(AR) )}. Thus, \var{ABS(AR)} is smaller than \const{SFMIN}. We can use this information and the assumption on \var{ABS(UR)} to conclude that \var{AR * AR + AI * AI} is smaller than or equal to \const{SFMIN * SFMIN}. Using the expression for \var{UI}, we also conclude that \var{ABS(UI)} is smaller than or equal to \const{(SFMIN * SFMIN) / ABS(AI)}. This means the absolute value \var{UI} is also small.
\end{remark}

\section{Numerical Analysis}
\newcommand{\CC}{\mathbb{C}}
\newcommand{\NN}{\mathbb{N}}
\newcommand{\RR}{\mathbb{R}}

This section gives a brief numerical analysis for the algorithm in the previous section. I analyze each of three cases separately. Let $x \in \CC^n$, $n \in \NN$, and $a \in \CC$ be given. The result of the algorithm is denoted by $y \in \CC^n$, and $a_r$ and $a_i$ are the real and imaginary parts of $a$, respectively. The algorithm produces a scaling factor $s$ which is a power of two.

Following \cite[§3.1]{Higham2002}, let
\begin{align*}
	\gamma_k := \frac{ku}{1-ku}\,, \quad k \in \NN\,,
\end{align*}
where $u$ is the unit roundoff. The relative error in the multiplication of $k$ numbers is bounded by $\gamma_k$, and this fact is used in the following analysis. Complexity is measured based on the number of additions and multiplications of floating-point numbers.

\subsection{Case 1: $a_i$ is zero}

In this case, the algorithm calls \func{RSCL} to scale $x$ by the reciprocal of $a_r$. The output is $y_i = (s/a_r)((1/s)x_i)$, $i=1,\ldots,n$. The factor $s$ is $1$ when $a_r$ is in the safe range. Therefore, the relative error in the imaginary and real parts of $y_i$ is bounded by $\gamma_2$. The complexity is either $O(2n)$, if $a_r$ is in the safe range, or $O(4n)$, in other cases.

\subsection{Case 2: $a_i$ is nonzero and $a_r$ is zero}

In this case, the output is $y_i = (0-(s/a_r)i)((1/s)x_i)$, $i=1,\ldots,n$. Analogously to the first case, the relative error in the imaginary and real parts of $y_i$ is bounded by $\gamma_2$. The complexity is either $O(2n)$, if $a_r$ is in the safe range, or $O(4n)$, in other cases. Notice that the algorithm for complex multiplications may not be optimized enough to avoid multiplications by zero. In that case, the complexity would be $O(6n)$. Moreover, such an algorithm for complex multiplication may generate NaNs in $y_i$ if $x_i$ is infinite.

\subsection{Case 3: $a_i$ and $a_r$ are both nonzero}

In this case, let $u_r, u_i \in \RR$ be defined as follows:
\begin{align*}
  u_r &= a_r + a_i(a_i/a_r), \\
  u_i &= a_i + a_r(a_r/a_i).
\end{align*}
In this case, the output is $y_i = b((1/s)x_i)$, where $b = (s/u_r) - (s/u_i)i$. The factor $s = 1$ when $u_r$ and $u_i$ are in the safe range. The relative error of $y_i$ is bounded by $\sqrt{2}\gamma_6$ (see \cite[§3.6]{Higham2002}). The complexity is either $O(6n)$, if $a_r$ is in the safe range, or $O(8n)$, in other cases.

Notice that the relative error bound $\sqrt{2}\gamma_6$ is ~50\% higher than the bound $\sqrt{2}\gamma_4$ for a straight-forward complex division (see \cite[§3.6]{Higham2002}). I believe this is acceptable provided that it is the same increase one can observe when analyzing \func{\{S,D,CS,ZD\}RSCL}. We lose accuracy to gain performance.

\section{Issues solved by the proposed algorithm}

\texttt{GETF2} is a LAPACK routine that computes the LU factorization of a matrix $A$ using partial pivoting with row interchanges. At every step $j$ of \texttt{GETF2}, the entries $A_{kj}$, $k > j$ are divided by $A_{jj}$, the chosen pivot. In the complex case, the routine needs to compute complex divisions to update these entries. The routines \{\texttt{C,Z}\}\texttt{GETF2} from LAPACK 3.11.0 may generate exceptions that are solved by the algorithm proposed in the previous section. See the examples below.

\begin{enumerate}
    \item Let
\[
A = \begin{bmatrix}
M+Mi & M \\
M & 0
\end{bmatrix},
\]
where $i=\sqrt{-1}$ and $M = 2^{127}$, which is approximately half of the overflow limit (OV) for single precision. Note that $A$ is very well conditioned, with $\|A\|_2 \|A^{-1}\|_2 \approx 3.7$. Note as well that $A$ is very close to overflow since $M$ is approximately OV/2. After running \texttt{CGETF2} from LAPACK 3.11.0 (gfortran 9.4.0 on \texttt{x86\_64-linux-gnu}) we obtain the factorization
\[
L = \begin{bmatrix}
1 & 0 \\
0 & 1
\end{bmatrix} \quad\text{ and }\quad
U = \begin{bmatrix}
M+Mi & M \\
0 & 0
\end{bmatrix},
\]
which raises an exception (\texttt{INFO = 2}) since $U$ is singular.
These $L$ and $U$ factors have a large backward error and are not correct.
The exact decomposition is
\[
L = \begin{bmatrix}
1 & 0 \\
0.5-0.5i & 1
\end{bmatrix} \quad\text{ and }\quad
U = \begin{bmatrix}
M+Mi & M \\
0 & -M(0.5-0.5i)
\end{bmatrix}.
\]
In this case, the main problem is the implementation of the complex division in the compiler~\footnote{\url{https://gcc.gnu.org/bugzilla/show_bug.cgi?id=107753}}.

\item Let $A$ be the ill-conditioned matrix
\[
A = \begin{bmatrix}
b+i & b \\
b & b
\end{bmatrix},
\]
where $b = 2^{75}$. Notice that $b$ is far below the overflow threshold, but $1/b^2$ underflows to 0. After running \texttt{CGETF2} from LAPACK 3.11.0 we obtain the factorization
\[
L = \begin{bmatrix}
1 & 0 \\
1 & 1
\end{bmatrix} \quad\text{ and }\quad
U = \begin{bmatrix}
b+i & b \\
0 & 0
\end{bmatrix},
\]
which again raises an exception (\texttt{INFO = 2}) since $U$ is singular. The exact decomposition is
\[
L = \begin{bmatrix}
1 & 0 \\
x-yi & 1
\end{bmatrix} \quad\text{ and }\quad
U = \begin{bmatrix}
b+i & b \\
0 & y+xi
\end{bmatrix},
\]
where $x = b^2/(b^2+1)$ and $y = b/(b^2+1)$. In single precision, $x$ and $y$ would round to $1$ and $1/b$ respectively, so that a high-accuracy answer would be
\[
L = \begin{bmatrix}
1 & 0 \\
1-i/b & 1
\end{bmatrix} \quad\text{ and }\quad
U = \begin{bmatrix}
b+i & b \\
0 & 1/b+i
\end{bmatrix}.
\]
In this case, the issue is $fl(fl(1/(b+i))b) = fl((1/b)b) = 1 \neq 1 - i/b = fl(b/(b+i))$.
Therefore, the result would be the same even if one uses \texttt{CLADIV} to compute $1/(b+i)$.

\end{enumerate}

\section{Conclusions}

This document describes an algorithm to scale a complex vector by the reciprocal of a complex constant. The algorithm has complexity O(6n) if scaling is not needed and O(8n) if scaling is needed, and does not perform any complex division. Although there is no performance comparisons in the text, we note that a naive approach of dividing two complex numbers would use at least 13 operations, including two floating-point divisions. The relative error in each entry of the computed vector is approximately 50\% higher than what is obtained with a straight-forward complex division. This is the same increase one observes when scaling by the reciprocal of a real number. Moreover, the new algorithm solves at least two issues in LAPACK 3.11.0 as detailed in this document.

\bibliographystyle{plain}
\bibliography{main}

\end{document}